\documentclass[submission]{FPSAC2027}

\usepackage{tikz, graphicx, float}
\usetikzlibrary{arrows,arrows.meta}
\definecolor{primary}{HTML}{2166AC}
\definecolor{secondary}{HTML}{2E7D32}
\definecolor{finalcomparison}{HTML}{C62828}
\definecolor{projectorblue}{RGB}{28,88,138}

\newcommand{\Z}{\mathbb{Z}}

\newcommand{\C}{\mathbb{C}}

\DeclareMathOperator{\height}{ht}

\DeclareMathOperator{\len}{len}
\DeclareMathOperator{\area}{area}
\DeclareMathOperator{\dinv}{dinv}
\DeclareMathOperator{\spec}{Spec}

\newcommand{\rot}{\geq_{\mathrm{rot}}}
\newcommand{\rotgeq}{\geq_{\mathrm{rot}}}
\newcommand{\rotleq}{\leq_{\mathrm{rot}}}

\newcommand{\uv}{\underline{v}}
\newcommand{\uov}{\underline{0v}}
\newcommand{\uz}{\underline{0}}
\newcommand{\Y}{\C[\mathbf{y}]}

\newcommand{\fp}{\mathfrak{p}}

\newcommand{\Hilb}{\mathrm{Hilb}}
\newcommand{\HHH}{\mathrm{HHH}}

\newcommand{\FT}{\mathrm{FT}}
\newcommand{\scd}{\mathrm{sc}}

\newtheorem{thm}{Theorem}
\newtheorem{lem}[thm]{Lemma}
\newtheorem{cor}[thm]{Corollary}
\newtheorem{defn}[thm]{Definition}
\newtheorem{example}[thm]{Example}
\newtheorem{remark}[thm]{Remark}
\newtheorem{prop}[thm]{Proposition}
\newtheorem{conj}[thm]{Conjecture}

\usepackage{lipsum}

\title[A filtered basis for Haiman ideals]{A filtered basis for Haiman ideals}

\author[Joshua P. Turner]{Joshua P. Turner \thanks{\href{mailto:jturner1295@gmail.com}{jturner1295@gmail.com}. This work was supported in part by the Pacific Institute for the Mathematical Sciences.}\addressmark{1}}

\address{\addressmark{1} Department of Mathematics, University of British Columbia}

\received{\today}

\abstract{
Given the polynomial ring $R=\C[x_1, \dots, x_n, y_1, \dots, y_n]$, we consider the Haiman ideal $J(d)$ of polynomials that vanish to order $d$ on the diagonals. We define a rotational order on monomials and use it to construct a filtration on $J(d)$ and a bigraded, homogeneous, free $\C[\mathbf{y}]$-basis compatible with this filtration. This gives a new independent proof that $J(d)$ is free over $y$'s. We also discuss the connections between this filtration and triply-graded link homology.
}


\usepackage[backend=bibtex]{biblatex}
\begin{document}

\maketitle

\section{Introduction}

Given the polynomial ring in two sets of variables $R_n = \C[x_1, \dots, x_n, y_1, \dots, y_n]$, we can consider the ideal $J_n(d) \subseteq R_n$, defined as
$$J_n(d) = \bigcap_{i<j} (x_i-x_j, y_i-y_j)^d.$$

We will suppress $n$ and $d$ when there is no ambiguity. This is the ideal of polynomials that vanish to order $d$ on each diagonal $x_i=x_j, y_i=y_j$. It has a natural bigrading, where the $x$'s and $y$'s have degree $q$ and $t$ respectively. This ideal was shown to be free over $\C[\mathbf{y}] = \C[y_1, \dots, y_n]$ by Haiman in his proof of the $n!$ and the Macdonald positivity conjectures \cite{Haiman}. His approach was to construct a filtration and a free $\C[\mathbf{y}]$-module basis for the larger polygraph ring $R(n,dn)$, and then show that $J(d)$ is isomorphic to the sign component of $R(n,dn)$. However this construction does not directly give a basis or bigraded Hilbert series for $J(d)$ itself.

Remarkably, we can recover the Hilbert series for $J(d)$ from computations of link homology. In \cite{Hog17}, Hogancamp uses a recursive process to compute the triply-graded link homology (abbreviated as $\HHH$) of the $d$'th power of a full twist as a power series in $q,t,$ and $a$. In \cite{GH}, Gorsky and Hogancamp identify the $a=0$ part of $\HHH$ of this braid with the quotient of $J(d)$ by the action of $y$'s. Putting these together, we can conclude that $J(d)$ is free over $y$'s with Hilbert series

\begin{equation}
\label{eqn: hilb series}
\Hilb_{q,t}\left(J_n(d)\right) = \frac{1}{(1-t)^n} \sum_{e \in \Z_{\geq 0}^n} q^{\area(e)}t^{\dinv^d(e)}.
\end{equation}

Our main result is the construction of a filtration on $J(d)$ and a compatible free $\Y$-basis that parallels this link homology computation.

\begin{thm}
    \label{thm: intro basis}
    There exists an order $\rot$ on nonnegative integer sequences, a filtration $N_v$ indexed by $v \in [0,d]^n$, and a bihomogeneous basis $p_e$ indexed by $e \in \Z_{\geq 0}^n$ with the following properties.
    \begin{enumerate}
        \item $N_{0^n} = J(d)$,
        \item $\deg(p_e) = (\area(e), \dinv^d(e))$,
        \item $\displaystyle p_e = c(\mathbf{y}) x^e + \sum_{f\, >_{\mathrm{rot}} \,e} c_f(\mathbf{y})x^f$,
        and 
        \item $N_v$ is a free $\Y$-module with basis $\{p_e \mid e \rot v\}$.
    \end{enumerate}
\end{thm}

This gives a new independent proof that $J(d)$ is free over $y$'s, and recovers the Hilbert series known from link homology.

\begin{cor}
    \label{cor: J free}
    For all $n,d$, the ideal $J(d)$ is a free $\Y$-module with bigraded Hilbert series
    $$\Hilb_{q,t}\left(J(d)\right) = \frac{1}{(1-t)^n} \sum_{e \in \Z_{\geq 0}^n} q^{\area(e)}t^{\dinv^d(e)}.$$
\end{cor}

In Section \ref{sec: rotation}, we define the order $\rot$ and the statistics on integer sequences $e$ that give the degrees of the basis elements $p_e$. In Section \ref{sec: filtration}, we define the filtration $N_v$, and in Theorem \ref{thm: filtration recursions}, we lay out the recursive properties of the associated graded modules $M_{\uv}$ that are required to show Theorem \ref{thm: intro basis}. Finally in Section \ref{sec: link homology}, we discuss in more detail the connections between our filtration and triply-graded link homology.

\section{Integer sequences and rotation}
\label{sec: rotation}

Given an integer sequence $e \in \Z_{\geq 0}^n$, first we define the statistics $\area(e)$ and $\dinv^d(e)$, which will give us the $q$ and $t$ degree of our basis elements $p_e$. For the sake of making some inductive steps clearer, we write integer sequences with decreasing indices, so $e = (e_n, \dots, e_1)$.

\begin{defn}
    For all $e \in \Z_{\geq 0}^n$, let
    $$\area(e) = |e| = \sum_{i \leq n} e_i,$$
    and let
    $$\dinv^d(e) = \sum_{i > j} \scd^d(e_i - e_j),$$
    where
    $$\scd^d(r) = 
    \begin{cases}
    d + r \; & \mathrm{if} \; -d \leq r \leq 0\\
    d + 1 - r \; & \mathrm{if} \; 1 \leq r \leq d\\
    0 \; & \mathrm{otherwise}
    \end{cases}$$
\end{defn}

The area and dinv statistics will give the $x$ and $y$ degrees of our basis elements $p_e$.

\begin{example}
    Let $d=1$. In this case, $\dinv = \dinv^1$ counts the number of pairs in $e$ that are either the same or decrease by 1. Consider the sequences $001$, $010$, and $0125410$. We have 
    $$\area(001) = 1, \quad \dinv(001) = 1,$$
    $$\area(010) = 1, \quad \dinv(010) = 2,$$
    and
    $$\area(0125410) = 13, \quad \dinv(0125410) = 6.$$ The dinv pairs for these examples are shown in Figure \ref{fig:dinv ex} below.
\end{example}

\begin{figure}[H]
    \centering
\scalebox{0.7}{
\begin{tikzpicture}[x=1.15cm,y=1.05cm,
  digit/.style={font=\Large,inner sep=3pt},
  equal/.style={draw=primary,line width=1.05pt},
  drop/.style={draw=secondary,line width=1.05pt}]

  \node[digit] (e1) at (-1,0) {$0$};
  \node[digit] (e2) at (0,0) {$0$};
  \node[digit] (e3) at (1,0) {$1$};
  \draw[equal] (e1.north) .. controls +(0,0.5) and +(0,0.5) ..
    node[above=3pt,font=\small,text=primary] {equal} (e2.north);

  \node[digit] (f1) at (4.3,0) {$0$};
  \node[digit] (f2) at (5.3,0) {$1$};
  \node[digit] (f3) at (6.3,0) {$0$};
  \draw[equal] (f1.north) .. controls +(0,0.7) and +(0,0.7) ..
    node[above=3pt,font=\small,text=primary] {equal} (f3.north);
  \draw[drop] (f2.south) .. controls +(0,-0.5) and +(0,-0.5) ..
    node[below=3pt,font=\small,text=secondary]
      {drop by $1$} (f3.south);

  \begin{scope}[shift={(8.6,0)}]
    \foreach \i/\val in {1/0,2/1,3/2,4/5,5/4,6/1,7/0}{
      \node[digit] (g\i) at ({\i-1},0) {$\val$};
    }

    \draw[equal] (g1.north)
      .. controls +(0,1.25) and +(0,1.25) .. (g7.north);
    \draw[equal] (g2.north)
      .. controls +(0,0.8) and +(0,0.8) .. (g6.north);

    \draw[drop] (g2.south)
      .. controls +(0,-1.55) and +(0,-1.55) .. (g7.south);
    \draw[drop] (g3.south)
      .. controls +(0,-1.0) and +(0,-1.0) .. (g6.south);
    \draw[drop] (g4.south)
      .. controls +(0,-0.45) and +(0,-0.45) .. (g5.south);
    \draw[drop] (g6.south)
      .. controls +(0,-0.45) and +(0,-0.45) .. (g7.south);

    \node[font=\large] at (3,-1.8)
      {$\operatorname{dinv}(0125410)=6$};
  \end{scope}

\end{tikzpicture}
}
\caption{Some dinv examples}
\label{fig:dinv ex}
\end{figure}
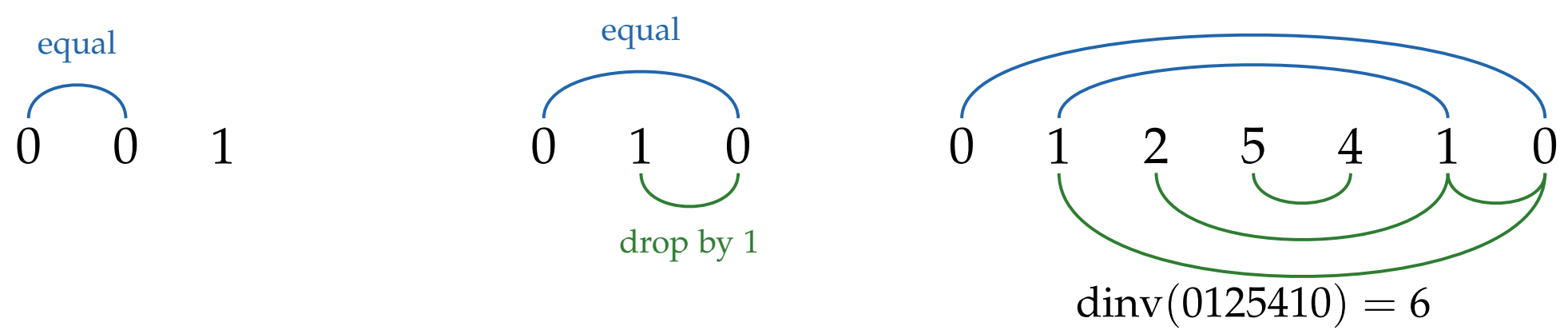

\begin{defn}
    Given an integer sequence $e \in \Z_{\geq 0}^n$, if $e_n > 0$, we define the rotation operation 
    $$\rho(e) := (e_{n-1}, \dots, e_1, e_n -1).$$ 
\end{defn}
This is the same rotation operation that appears in \cite{Hog17} and very similar to the one in \cite{Haiman}. 
\begin{remark}
    The rotation operation $\rho$ lowers area by 1 and preserves dinv. Indeed,
    $$\area(\rho(e)) = \area(e) - 1, \; \textrm{and} \; \dinv^d(\rho(e)) = \dinv^d(e)$$
    for all $e$ with $e_n > 0$.
\end{remark}
\begin{defn}
    We define an order on integer sequences $e \in \Z_{\geq 0}^n$ by induction as follows.
    \begin{enumerate}
        \item If $e_n, f_n >0$, then
        $$e \rot f \iff \rho(e) \rot \rho(f).$$
        \item If $e_n = f_n = 0$, and $e = 0r$, $f = 0s$
        $$e \rot f \iff r \rot s.$$
        \item If $e_n > 0$ and $f_n = 0$, then
        $$e \rot f.$$
    \end{enumerate}
\end{defn}

It is clear that this is well defined inductively, as at every step we either perform a rotation (i.e. apply $\rho$), which decreases $\area(e)$ and $\area(f)$, or we perform a `closing up' operation (i.e. eliminate a leading $0$), which decreases $\len(e)$ and $ \len(f)$. To compare any two sequences, we apply these steps until exactly one of the two sequences starts with $0$, or both are empty. Also note that although this order is not lexicographic, it is true that if $e_i \geq f_i$ for all $i$, then $e \rot f$.

Alternatively, we can encode the sequence of rotation and closing up steps that take any sequence $e$ to the empty sequence, and we can compare two sequences by comparing their binary sequences.

\begin{defn}
\label{def: binary rot order}
    For any sequence in  $\Z_{\geq 0}^n$, define the binary sequence $r(e)$ inductively as follows:
    \begin{enumerate}
        \item $r(\varnothing) = \varnothing$;
        \item $r(0e) = 0r(e)$;
        \item If $e_n >0$, $r(e) = 1 r(\rho(e))$
    \end{enumerate}
\end{defn}

\begin{prop}
\label{prop: binary rot order}
    For any $e, f \in \Z_{\geq 0}^n$, 
    $$e \rot f \iff r(e) \geq_{\mathrm{lex}} r(f).$$
\end{prop}

\begin{example}
    \label{rot example}
    Below, we compare the sequences $001, 010$ and $011, 020$ using both methods. In general, sequences that have larger area and are `flatter' tend to be larger.
\end{example}

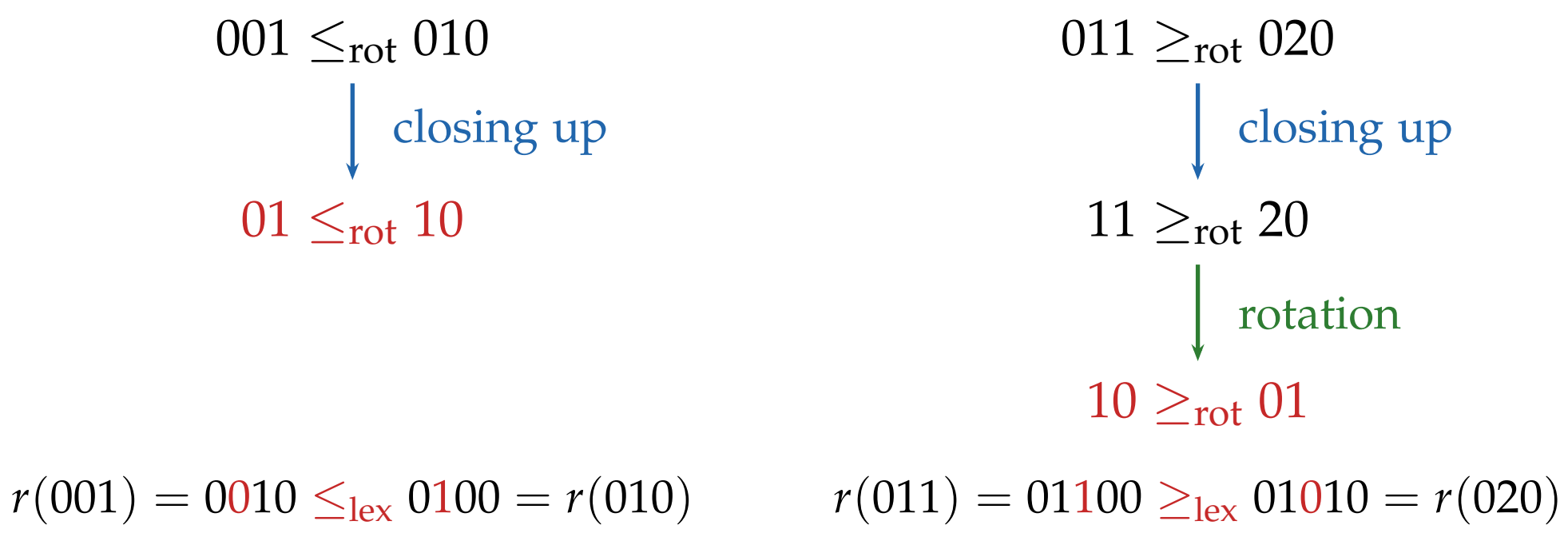
\begin{figure}[H]
    \centering
\begin{tikzpicture}[x=1cm,y=1cm,
  reduction/.style={-{Stealth[length=4pt]},line width=1pt,primary},
  rotation/.style={reduction,secondary}]
  \node at (-2,2.8) {$001\rotleq010$};
  \node at (5,2.8) {$011\rotgeq020$};
  \node[text=finalcomparison] at (-2,1.3) {$01\rotleq10$};
  \node at (5,1.3) {$11\rotgeq20$};
  \node[text=finalcomparison] at (5,-0.2) {$10\rotgeq01$};
  \foreach \x in {-2,5} {
    \draw[reduction] (\x,2.45) -- (\x,1.65)
      node[midway,right=5pt,font=\small] {closing up};
  }
  \draw[rotation] (5,0.95) -- (5,0.15)
    node[midway,right=5pt,font=\small] {rotation};
    
  \node[font=\small] at (-2,-1)
    {$r(001)=0\textcolor{finalcomparison}{0}10
      \mathrel{\textcolor{finalcomparison}{\leq_{\mathrm{lex}}}}
      0\textcolor{finalcomparison}{1}00=r(010)$};
  \node[font=\small] at (5,-1)
    {$r(011)=01\textcolor{finalcomparison}{1}00
      \mathrel{\textcolor{finalcomparison}{\geq_{\mathrm{lex}}}}
      01\textcolor{finalcomparison}{0}10=r(020)$};
\end{tikzpicture}
\caption{Comparing sequences with $\rot$}
\label{fig: comparing with rot}
\end{figure}

\section{Filtration and basis}
\label{sec: filtration}

Our filtration for $J(d)$ will be indexed by bounded sequences $v \in [0,d]^n$, which can be compared using $\rot$ just as unbounded sequences.

\begin{defn}
\label{defn: Fv and Mv}
    For any $v \in [0,d]^n$, define the monomial ideal
    $$F_v = (x^e \mid e \rot v) \subseteq R,$$ and let
    $$N_v = F_v \cap J(d).$$
    Here $x^e = \prod_{i \leq n} x_i^{e_i}$.
\end{defn}

 Essentially, we have an order on monomials in $x$ induced by $\rot$, and our filtration consists of certain level sets of this order inside $J(d)$.

\begin{remark}[Warning!]
It is important to note that the $\rot$ order on monomials in $x$ is \textbf{not} a true monomial order. Indeed, Example \ref{rot example} demonstrates that it is not multiplicative, as $$x^{001} \rotleq x^{010}, \; \textrm{ but } \, x_2 (x^{001}) \rot x_2 (x^{010}).$$

Thus our filtered basis is not a Groebner basis.
\end{remark}

It's clear from the definition that if $w \rot v$, then $N_w \subseteq N_v$, and that $N_{0^n} = J(d)$, so the modules $N_v$ indeed form a filtration of $J(d)$. Since $\rot$ is a total order on the finite set of bounded sequences, every $v \in [0,d]^n$ has a unique successor $v^+$ other than the unique maximal element $v = d^n$. 

\begin{defn}
    Given $v \in [0,d]^n$, $v\neq d^n$, let $v^+$ be the unique successor of $v$. Define the associated graded module
    $$M_{\uv} := N_v/N_{v^+}.$$
    In the case where $v = d^n$, for simplicity we define
    $$M_{d^n} = N_{d^n}.$$
\end{defn}

Now we characterize which monomials in $x$ survive in these associated graded pieces. For $e \in \Z_{\geq 0}^n$ and $v \in [0,d]^n$, we say that $e$ \textit{covers} $v$ if for all $i \leq n$, if $v_i<d$, then $e_i = v_i$, and if $v_i = d$, then $e_i \geq v_i$.

\begin{prop}
    \label{prop: e covers v}
    For an unbounded sequence $e \in \Z_{\geq 0}^n$ and a bounded sequence $v \neq d^n$, $e$ covers $v$ if and only if $e \rot v$ and $e <_{\mathrm{rot}} v^+$.
\end{prop}

In other words, the monomials $x^e$ that are in $F_v$ but not $F_{v+}$ are exactly those where $e$ covers $v$. We can go further to find a different characterization for $M_{\uv}$.

\begin{prop}
    \label{prop: Mw = Xv cap Mv}
    Let $X_v$ be the monomial ideal
    $$X_v = (x_i^{v_i+1} \mid v_i < d).$$
    If $w = v^+$ is the successor of $v$, then $F_w = X_v \cap F_v$ and $N_w = X_v \cap N_v$.
\end{prop}

\begin{cor}
\label{cor: Mv + Xv}
    The associated graded module $M_{\uv}$ is isomorphic to $\left(N_v + X_v \right)/X_v.$
\end{cor}

Given this isomorphism, we can interpret $M_{\uv}$ as a module over $R_n/X_v$. The motivation behind the underline notation is that $\uv$ is shorthand for the sequence 
$$\uv = \begin{cases}
         \underline{v_i} & \text{if } v_i < d \\
        d  & \text{if } v_i = d
        \end{cases}$$
in the alphabet $\{\underline{0}, \underline{1}, \dots, \underline{d-1}, d\}$. The monomials in $x$ in $N_v$ that survive mod $X_v$ are precisely the monomials of the form $x^e$ where $e \rot v$ and where any underlines in $v$ impose the extra condition that $e_i = \underline{v_i}$. This is equivalent to the condition that $e$ covers $v$.

\begin{example}
    Let $n=3$ and $d=1$. Consider $v = 001$. In this case, we have that $$F_{001} = (x^e \mid e \rot 001) = (x_1, x_2, x_3)$$ and $$X_{001} = (x_2, x_3).$$
    
    Then
    $$M_{\uv} = M_{\underline{00}1} = N_{001}/N_{010} \cong \left(N_{001} + (x_2, x_3)\right)/(x_2, x_3)$$

    Thus we can think of $M_{\underline{00}1}$ as the space of polynomials in $(x_1) \subseteq R/(x_2,x_3) = \C[x_1, y_1, y_2, y_3]$ that are in the image of $J$. See Section \ref{sec: n=3} for a more detailed characterization of $M_{\underline{00}1}$.
\end{example}

\subsection{Filtration recursions}

Our goal is to show by induction that the associated graded modules $M_{\uv}$ are free over $y$'s, and to construct a free $\Y$-module basis of the form $\{p_e \mid e \textrm{ covers } v\}$. We do this by induction on $v$, with the following 3 inductive steps.

\begin{thm}
    \label{thm: filtration recursions}
    Let $v \in [0,d]^{n-1}.$ Let $\pi$ be the cyclic permutation $\pi = (n-1 \ \dots \,  1 \quad n)$ that acts on $R$ by permuting $x$'s and $y's$ simultaneously, and let 
    $$\Delta_v = \prod_{i < n} (y_i - y_n)^{d-v_i}.$$
    We have the following recursions for the associated graded modules $M_{\uv}$.
    \begin{enumerate}
        \item (Rotation) If $1 \leq k < d$, then
        $$M_{\underline{kv}} = x_n \pi(M_{\underline{v(k-1)}})$$
        \item (Splitting) There is a short exact sequence of $R$-modules
        $$0 \rightarrow x_n\pi\left( M_{\uv d}\right) \rightarrow M_{d \uv} \rightarrow x_n \pi \left( M_{\underline{v(d-1)}}\right) \rightarrow 0$$
        \item (Closing Up) We have the following isomorphism of $R_n$ modules.
    \begin{equation}
        M_{\uov} \cong \Delta_v M_{\uv} \otimes \C[y_n]
    \end{equation}
    Here $M_{\uov}$ is an $R_n$-module (with a trivial action of $x_n$) and $M_{\uv}$ is an $R_{n-1}$-module.
    \end{enumerate}
\end{thm}
\begin{proof}[Proof sketch of Steps (1) and (2)]
    Since $J(d)$ itself is invariant under the action of the permutation $\pi$, Step (1) essentially reduces to verifying that the allowed monomials are the same on both sides. Note that for any $e$, $x^e = x_n \pi (x^{\rho(e)}),$ so this step is the operation that rotates every monomial in $M_{\underline{kv}}$ simultaneously. It is sufficient to check that $e$ covers $kv$ if and only if $\rho(e)$ covers $\rho(kv) = v(k-1)$.

    Similarly, Step (2) reduces to checking that $e$ covers $d\underline{v}$ if and only if $\rho(e)$ covers $\underline{v}d$ or $\underline{v(d-1)}$. The splitting occurs from separating the case where $e_n = d$ from the case where $e_n > d$, which lead to two different truncated sequences $v$ after rotation.
\end{proof}

The `closing up' Step (3) tells us that $M_{\uov}$ is isomorphic to $M_{\uv}$ up to a degree shift in $y$'s. This means that given a basis of the form $\{q_{e} \mid e \textrm{ covers } \uv \}$, the polynomials $q_{0e} = \Delta_v q_e$ form a basis for $M_{\uov}$. Therefore, Step (3) is equivalent to the statement that for all $e$ that cover $v$, there exists a lift $p_{0e} \in N_{0v}$ of $q_{0e}$ such that 
$$p_{0e} = q_{0e} = \Delta_v q_e \textrm{ (mod }X_{0v}).$$

Finding a concrete combinatorial way to perform this lift has proven difficult, so instead we resort to using non-constructive localization arguments. See Section \ref{sec: closing up} for a detailed proof sketch of Step (3).

Once we have these filtration recursions, our main result Theorem \ref{thm: intro basis} follows by induction. Every step in Theorem \ref{thm: filtration recursions} either reduces $\area(v)$ or reduces $n$. For steps (1) and (3), if we assume by induction that the right hand side is free over $y$'s with the required basis, we can check that the left hand side will also be free with the desired basis $p_e$. For step (2), if we assume by induction that both $M_{\uv d}$ and $M_{\underline{v(d-1)}}$ are free over $y$'s with the desired basis, then the short exact sequence splits over $\Y$, and as $\Y$-modules,
\begin{equation}
    M_{d \uv} \cong x_n \pi \left(M_{\uv d} \oplus M_{\underline{v(d-1)}}\right).
\end{equation}
Lastly, we verify that the degree shifts from these filtration recursion steps match the degree shifts in $q$ and $t$ that come from rotation and closing up of integer sequences. Rotation simply shifts the $q$ degree by $1$, and closing up shifts the $t$ degree by 
$$\dinv^d(0v) - \dinv^d(v) = \deg_{t}(\Delta_v).$$

\subsection{Filtration for $n=3$}
\label{sec: n=3}
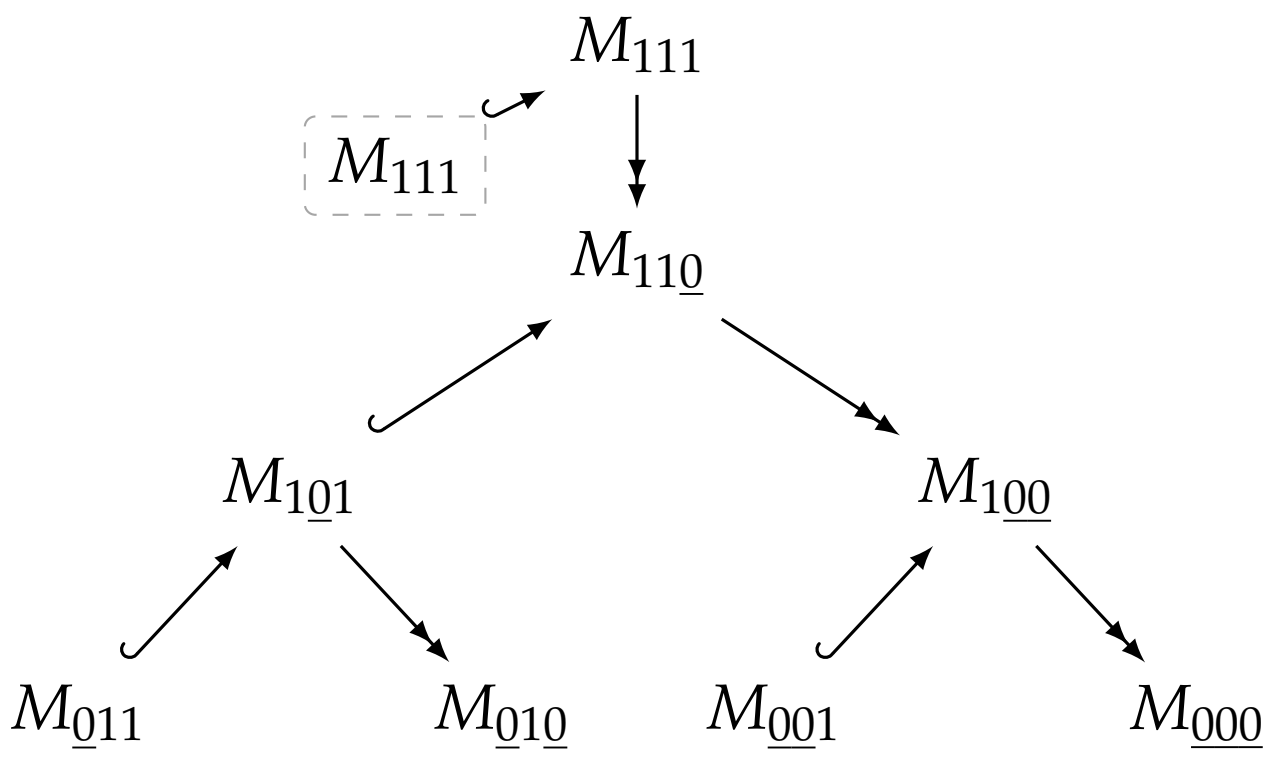
\begin{figure}[H]
    \centering
\begin{tikzpicture}[
  x=1cm,y=1cm,>=latex,
  every node/.style={inner sep=4pt},
  inc/.style={right hook->,semithick},
  surj/.style={->>,semithick},
  leaf/.style={draw=black!35,fill=black!3,
               rounded corners=2pt,inner sep=6pt}
]
  \node (root) at (0.5,0) {$M_{111}$};
  \node[draw=black!35,dashed,rounded corners=2pt]
        (repeat) at (-1.1,-.8) {$M_{111}$};
  \node (a) at (0.5,-1.45) {$M_{11\uz}$};

  \node (b) at (-1.8,-2.95) {$M_{1\uz1}$};
  \node (c) at (2.8,-2.95) {$M_{1\uz\uz}$};

  \node (d) at (-3.2,-4.45) {$M_{\uz11}$};
  \node (e) at (-0.4,-4.45) {$M_{\uz1\uz}$};
  \node (f) at (1.4,-4.45) {$M_{\uz\uz1}$};
  \node (g) at (4.2,-4.45) {$M_{\uz\uz\uz}$};

  \draw[inc] (repeat) -- (root);
  \draw[surj] (root) -- (a);
  \draw[inc] (b) -- (a);
  \draw[surj] (a) -- (c);
  \draw[inc] (d) -- (b);
  \draw[surj] (b) -- (e);
  \draw[inc] (f) -- (c);
  \draw[surj] (c) -- (g);
\end{tikzpicture}
\caption{Filtration for $n=3, \, d=1$}
\label{fig: n=3 filtration}
\end{figure}

    Figure \ref{fig: n=3 filtration} below is a diagram of the filtration for $J = J_3(1)$. Each branching node represents the short exact sequence from Step (2), with rotations and grading shifts suppressed.

In order to find the basis for $J_3(1)$, we start with $M_{111}$, which is equal to $J$ up to a grading shift,
$$M_{111} = (x_1 x_2x_3) \cap J = x_1x_2x_3J.$$
We apply Step 2 to get the short exact sequence
$$0 \rightarrow x_n\pi\left( M_{111}\right) \rightarrow M_{111} \rightarrow x_n \pi \left( M_{11\underline{0}}\right) \rightarrow 0.$$
Since $M_{111}$ is invariant under $\pi$, we have that
$$x_3\pi\left(M_{11\underline{0}}\right) \cong M_{111}/(x_3)M_{111}.$$
In fact, $x_n$ acts freely on $M_{111}$, so we can easily recover the basis for $M_{111}$ from the basis for $M_{11\underline{0}}$. From here, we repeatedly apply Step (2) until we arrive at modules of the form $M_{\underline{0v}}$ with $v \in [0,1]^2$, represented by the leaves of the tree in Figure \ref{fig: n=3 filtration}.

At this point, we apply Step (3), which tells us for example that

$$M_{\underline{00}1} \cong (y_2 - y_3)M_{\underline{0}1} \otimes \C[y_3].$$

It's easy to check that $\{q_{0k} = x_1^k \mid k\geq 1\}$ forms a basis for $M_{\underline{0}1}$ (with $p_{0k} = x_1^{k-1}(x_1-x_2)$ being the lifted element in $J_2$). So Step (3) tells us that there exists some $p_{00k} \in N_{001}$ such that
$$p_{00k} = q_{00k} = (y_2-y_3)q_{0k} = (y_2-y_3)x_1^k \quad \textrm{ (mod }(x_2,x_3)),$$
but it does not tell us how to construct it. In this case we can check that
$$p_{001} = 
\begin{vmatrix}
1   & 1   & 1   \\
x_1 & x_2 & x_3 \\
y_1 & y_2 & y_3
\end{vmatrix},
$$
and
$$p_{00k} = x_1^{k-1} p_{001}$$

have the desired properties, and thus form a basis for $M_{\underline{00}1}$.

\subsection{Proof sketch of Closing Up}
\label{sec: closing up}

In this section, we provide an outline of the proof of Step 3 (Closing Up) of Theorem \ref{thm: filtration recursions}. We present it here again for convenience.

\begin{thm}[Closing Up]
    \label{thm: closing up}
    Let $v \in [0,d]^{n-1}$ and let $\Delta_v = \prod_{i < n} (y_i - y_n)^{d-v_i}.$
    Then we have that
    $$M_{\uov} \cong \Delta_v M_{\uv} \otimes \C[y_n]$$
    as $R_n$-modules.
\end{thm}

We start by constructing a natural injection
$$\phi_v : M_{\uov} \rightarrow \Delta_v M_{\uv} \otimes \C[y_n]$$ induced by specializing $x_n$ to be 0, and confirming that its image lies in $\Delta_v M_{\uv} \otimes \C[y_n]$. To show that $\phi_v$ is surjective, it would be sufficient to construct an explicit preimage of $\Delta_v p_e$ for every basis element $p_e$ of $M_{\uv}$. However, this has proven to be a difficult combinatorial problem, so we are forced to use localization arguments instead. We next show that $\phi_v$ is an isomorphism outside of codimension $3$ in $\spec(\Y)$.

\begin{defn}
    Define $Z \subseteq \spec(\Y)$ to be the union of all linear subspaces of the form $V(y_{i_1} - y_{j_1}, y_{i_2} - y_{j_2}, y_{i_3} - y_{j_3})$ with all 3 pairs required to be linearly independent. Let $U = \spec(\Y) \backslash Z$ be its complement.
\end{defn}
The open set $U \subseteq \spec(\Y)$ is the locus where the coordinates $y_1, \dots, y_n$ take on at least $n-2$ different values. This $U$ is the same set that Haiman labels $\hat{U}_3$ in \cite[Definition 4.4.1]{Haiman}. This $U$ is the natural choice of open set in $\spec(\Y)$ whose complement $Z$ has codimension 3 that plays nicely with the ideal $J$.
\begin{prop}
    \label{prop: local closing up}
    Theorem \ref{thm: closing up} holds locally on $U$. In particular, for all $v$, the injection $\phi_v$ becomes an isomorphism after restricting to $U$.
\end{prop}

To prove Proposition \ref{prop: local closing up}, we show that after localizing by a prime ideal $\fp \in U$, the basis lifting required to show that $\phi_v$ is surjective essentially reduces to the $n=3$ or $n=2$ case. Then we construct the required lifts in these cases explicitly. In Section \ref{sec: n=3}, we discuss how to do one of these lifts in the $n=3, d=1$ case.

Next, we argue that if $\phi_v$ is an isomorphism outside codimension 3, it is in fact an isomorphism globally using the following lemma, which is a consequence of the depth lemma for a short exact sequence. Here S2 and S3 are Serre's conditions for modules.

\begin{lem}
    \label{lem: reflexive SES}
    Suppose that $A$ is a regular integral domain and we have an exact sequence
    $$0 \rightarrow P \rightarrow M \rightarrow N \rightarrow 0$$
    of finitely generated, torsion-free $A$-modules. If
    \begin{enumerate}
        \item $P$ is free (or just S3),
        \item $M$ is S2, and
        \item $N_p$ is free over $A_{\fp}$ for every $\fp \in Spec(A)$ with $\height \fp \leq 2$,
    \end{enumerate}
    then $N$ is S2.
\end{lem}

We show that the short exact sequence

$$0 \rightarrow N_{(0v)^+} \rightarrow N_{0v} \rightarrow M_{\uov} \rightarrow 0$$
satisfies the conditions of Lemma \ref{lem: reflexive SES}. In particular, we assume that $N_{0w}$ and $M_{\uv}$ are free over $A = \Y$ by induction, and Proposition \ref{prop: local closing up} tells us that $(M_{\uov})_\fp$ is free for every $\fp$ of height at most 2. To make sure everything is finitely generated, we apply the lemma in each fixed $x$-degree.

Finally, we use Hartogs Lemma to conclude that since $M_{\uov}$ is S2 and $\phi_v$ is an isomorphism outside of codimension 2, it is in fact globally an isomorphism.
 
\section{Connections to link homology}
\label{sec: link homology}

By work of Gorsky and Hogancamp \cite{GH}, it is known that $\HHH^{a=0}(\FT_n^d) \cong J_n(d)/(y)J_n(d)$. In \cite{Hog17}, Hogancamp computed $\HHH$ of the braid $\beta = \FT_n^d$ using methods developed by Elias and Hogancamp in \cite{EH}. To do this computation, Hogancamp introduces projectors $K_n$ and auxiliary complexes $C_v$ indexed by $v \in [0,d]^n$ whose homology satisfies certain recursions that strongly parallel those in Theorem \ref{thm: filtration recursions}. Figure \ref{fig: Braid recursion} below shows the structure of these recursions in the case $n=3, d=1$.

\newcommand{\BraidStep}[3]{%
  \ifnum#1=0
    \foreach \xx in {-0.28,0,0.28}
      \draw[strand] (\xx,#2) -- (\xx,#3);
  \else
    \pgfmathsetmacro{\xl}{0.28*(#1-2)}
    \pgfmathsetmacro{\xr}{\xl+0.28}
    \pgfmathsetmacro{\yc}{0.55*(#2)+0.45*(#3)}
    \pgfmathsetmacro{\yd}{0.45*(#2)+0.55*(#3)}
    \ifnum#1=1 \def\xs{0.28}\else\def\xs{-0.28}\fi
    \draw[strand] (\xs,#2) -- (\xs,#3);
    \draw[strand] (\xr,#2)
      .. controls (\xr,\yc) and (\xl,\yd) .. (\xl,#3);
    \draw[strand,preaction={draw=white,line width=1.0pt}]
      (\xl,#2) .. controls (\xl,\yc) and (\xr,\yd) .. (\xr,#3);
  \fi
}

\newcommand{\BraidShuffle}[3]{%
  \pgfmathsetmacro{\yc}{0.65*(#2)+0.35*(#3)}
  \pgfmathsetmacro{\yd}{0.35*(#2)+0.65*(#3)}
  \ifnum#1=1
    \foreach \xx in {0,0.28} {
      \draw[strand] (\xx,#2)
        .. controls (\xx,\yc) and ({\xx-0.28},\yd)
        .. ({\xx-0.28},#3);
    }
    \draw[strand,preaction={draw=white,line width=1.0pt}]
      (-0.28,#2) .. controls (-0.28,\yc) and (0.28,\yd) .. (0.28,#3);
  \else
    \draw[strand] (0.28,#2)
      .. controls (0.28,\yc) and (-0.28,\yd) .. (-0.28,#3);
    \foreach \xx in {-0.28,0} {
      \draw[strand,preaction={draw=white,line width=1.0pt}]
        (\xx,#2) .. controls (\xx,\yc) and ({\xx+0.28},\yd)
        .. ({\xx+0.28},#3);
    }
  \fi
}

\newcommand{\BraidLabels}[3]{%
  \foreach \xx/\digit in {-0.28/#1,0/#2,0.28/#3} {
    \ifnum\digit=0
      \node[binary] at (\xx,0.35) {$\underline{0}$};
    \else
      \node[binary] at (\xx,0.35) {$1$};
    \fi
  }
}

\newcommand{\TwistBack}[2]{%
  \draw[ring] ({#1+#2},-0.09)
    arc[start angle=0,end angle=180,x radius=#2 cm,y radius=0.08cm];
}
\newcommand{\TwistFront}[2]{%
  \draw[ring,preaction={draw=white,line width=1.4pt}]
    ({#1-#2},-0.09)
    arc[start angle=180,end angle=360,x radius=#2 cm,y radius=0.08cm];
}

\newcommand{\MiniBraid}[1]{%
  \begin{tikzpicture}[
    x=1cm,y=1.2cm,baseline=-0.5ex,
    strand/.style={solid,draw=black,line width=0.45pt,line cap=round},
    ring/.style={solid,sharp corners,draw=black,line width=0.45pt},
    every node/.style={font=\fontsize{8}{9}\selectfont,
                       inner sep=1.1pt,outer sep=0pt},
    binary/.style={font=\fontsize{9}{10}\selectfont,anchor=base,inner sep=0pt},
    block/.style={solid,sharp corners,draw=black!75,fill=white,line width=0.4pt,
                  minimum height=0.28cm},
    projector/.style={block,draw=projectorblue,
                       text=projectorblue,fill=projectorblue!5,
                       font=\fontsize{6}{7}\selectfont,inner xsep=0.15pt}
  ]
    \useasboundingbox (-0.47,-0.53) rectangle (0.47,0.53);
    \def\ba{0}\def\bb{0}\def\nk{0}
    \ifnum#1=110 \def\ba{2}\def\bb{1}\def\nk{1}\fi
    \ifnum#1=101 \def\ba{1}\def\nk{1}\fi
    \ifnum#1=100 \def\ba{1}\def\bb{2}\def\nk{2}\fi
    \ifnum#1=011 \def\nk{1}\fi
    \ifnum#1=010 \def\ba{2}\def\nk{2}\fi
    \ifnum#1=001 \def\nk{2}\fi
    \ifnum#1=000 \def\nk{3}\fi

    \BraidLabels#1
    \ifnum\nk=0 \TwistBack{0}{0.39}\fi
    \ifnum\nk=1 \TwistBack{0.14}{0.20}\fi

    \ifnum\bb=0
      \BraidStep{\ba}{0.28}{0.06}
      \BraidStep{\ba}{-0.24}{-0.50}
    \else
      \BraidShuffle{\ba}{0.28}{0.06}
      \BraidShuffle{\bb}{-0.24}{-0.50}
    \fi
    \foreach \xx in {-0.28,0,0.28}
      \draw[strand,preaction={draw=white,line width=1.3pt}]
        (\xx,0.06) -- (\xx,-0.24);

    \ifcase\nk
      \TwistFront{0}{0.39}
    \or
      \node[projector,minimum width=0.22cm] at (-0.28,-0.09) {$K_1$};
      \TwistFront{0.14}{0.20}
    \or
      \node[projector,minimum width=0.50cm] at (-0.14,-0.09) {$K_2$};
    \or
      \node[projector,minimum width=0.76cm] at (0,-0.09) {$K_3$};
    \fi
  \end{tikzpicture}%
}
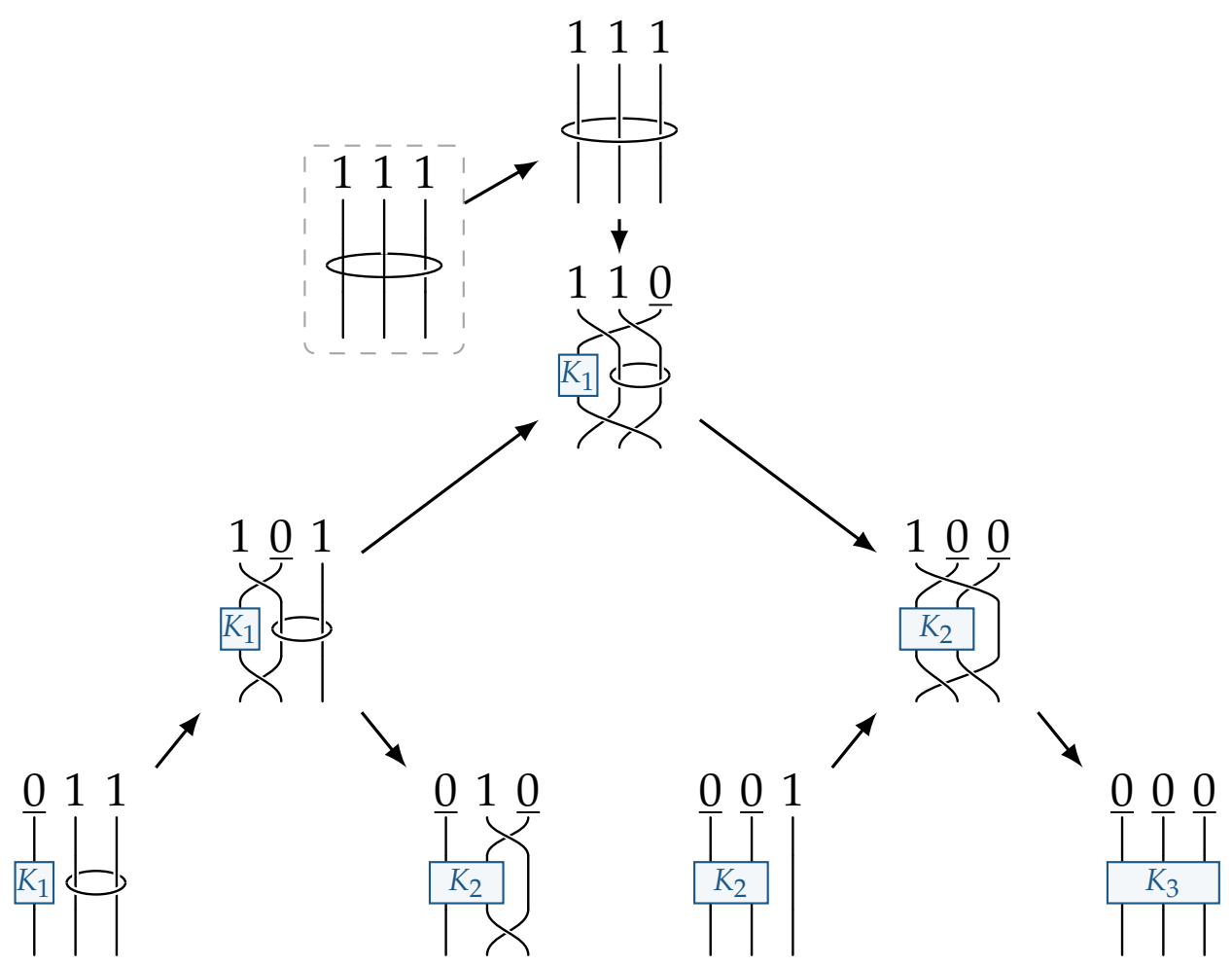
\begin{figure}[H]
    \centering
\begin{tikzpicture}[
  x=1cm,y=1.15cm,>=latex,
  every node/.style={inner sep=2pt},
  morphism/.style={->,semithick}
]
  \node (root) at (0.5,0) {\MiniBraid{111}};
  \node[draw=black!35,dashed,rounded corners=2pt]
        (repeat) at (-1.1,-.8) {\MiniBraid{111}};
  \node (a) at (0.5,-1.45) {\MiniBraid{110}};

  \node (b) at (-1.8,-2.95) {\MiniBraid{101}};
  \node (c) at (2.8,-2.95) {\MiniBraid{100}};

  \node (d) at (-3.2,-4.45) {\MiniBraid{011}};
  \node (e) at (-0.4,-4.45) {\MiniBraid{010}};
  \node (f) at (1.4,-4.45) {\MiniBraid{001}};
  \node (g) at (4.2,-4.45) {\MiniBraid{000}};

  \draw[morphism] (repeat) -- (root);
  \draw[morphism] (root) -- (a);
  \draw[morphism] (b) -- (a);
  \draw[morphism] (a) -- (c);
  \draw[morphism] (d) -- (b);
  \draw[morphism] (b) -- (e);
  \draw[morphism] (f) -- (c);
  \draw[morphism] (c) -- (g);

\end{tikzpicture}
    \caption{Braid recursion for $\FT_3$}
    \label{fig: Braid recursion}
\end{figure}

Let $m_v(q,t) = \Hilb_{q,t}(M_{\uv}/(\mathbf{y})M_{\uv})$, and let $l(v) = \deg_t(\Delta_v) = \dinv^d(0v)-\dinv^d(v)$. If we take the module recursions from Theorem \ref{thm: filtration recursions} and only keep track of the bigraded Hilbert series, we get the following recursions for $m_v$.
\begin{prop}
\label{prop: fv recursions}
Let $v \in [0,d]^n$. The series $m_v$ defined above satisfies the following.
\begin{enumerate}
    \item (Base case) $$m_{\varnothing} = 1$$
    \item (Rotation) If $1 \leq k < d$, then
    $$m_{kv} = q m_{v(k-1)}$$
    \item (Splitting) $$m_{dv} = qm_{vd} + qm_{v(d-1)}$$
    \item (Closing up) $$m_{0v} = t^{l(v)}m_v$$
\end{enumerate}    
\end{prop}

\begin{prop}
\label{prop: Cv hilbert series}
    Up to degree shifts, $m_v$ satisfies the same recursions as the series $f_v(q,t,0)$ in \cite[Theorem 1.3]{Hog17}.
\end{prop}

In particular this means that up to degree shifts, $\HHH^{a=0}(C_v)$ and $M_{\uv}/(\mathbf{y})M_{\uv}$ have the same bigraded Hilbert series. We expect that this actually holds on the level of bigraded $R$-modules.

\begin{conj}
    Up to degree shifts,
    $\HHH^{a=0}(C_v) \cong M_{\uv}/(\mathbf{y})M_{\uv}$
    as $R$-modules.
\end{conj}

We also hope to find a similar filtration for the generalized Haiman ideal
$$J(d_1, \dots, d_n) = \bigcap_{i<j} (x_i-x_j, y_i-y_j)^{\min(d_i,d_j)}.$$ It is conjectured that this ideal is free over $y$'s and should be related to $\HHH$ of certain Coxeter braids.

\acknowledgements{This work was supported in part by the Pacific Institute for the Mathematical Sciences. The author would like to thank Eugene Gorsky for helpful conversations.}

\nocite{*}
\printbibliography

\clearpage
\section*{AI disclosure}

\begin{itemize}
    \item Mathematical content: The statement and proof of Lemma \ref{lem: reflexive SES} were suggested by ChatGPT (GPT-5.6 Pro) during dialogue about completing the proof of Theorem \ref{thm: closing up}. The author verified and recreated the proof and revised its statement.
    \item Exposition: The author used ChatGPT (GPT-6 Astra) to assist in outlining and organizing this extended abstract, as well as in proofreading and final editing. However, all text was fully written by the author.
    \item Figures: All figures were generated with assistance from ChatGPT (GPT-6 Astra). They were generated from specifications provided by the author and were checked and edited directly by the author afterwards.
\end{itemize}
\end{document}